\documentclass[11pt]{amsart}
 \usepackage{graphicx}
 \usepackage{amssymb}

\newtheorem{theorem}{Theorem}[section]        
            
\newtheorem{corollary}[theorem]{Corollary}     
 
\newtheorem{problem}[theorem]{Open Problems} 
\newtheorem{definition}[theorem]{Definition}   
\newcommand{\nc}{\newcommand}
\nc{\Q}{\ensuremath{\mathbb{Q}}}
\nc{\R}{\ensuremath{\mathbb{R}}}
\nc{\M}{\ensuremath{\mathcal{M}}}
\nc{\F}{\ensuremath{\mathcal{F_{\text{acet}}}}}
\nc{\G}{\ensuremath{\mathcal{G}}}
\nc{\Or}{\ensuremath{\mathcal{O}}}
\nc{\N}{\ensuremath{\mathcal{N}}}
\nc{\Pa}{\ensuremath{\mathcal{P}}}
\nc{\La}{\ensuremath{\mathcal{L}}}
\nc{\Sa}{\ensuremath{\mathcal{S}}}
\nc{\C}{\ensuremath{\mathcal{C}}}
\nc{\Z}{\ensuremath{\mathbb{Z}}}
 \begin{document}       
 \title{Kalai  orientations on matroid polytopes}
 \author[Raul Cordovil]{Raul Cordovil
  \address{Departamento de Matem\'atica,\newline
Instituto Superior T\' ecnico,\newline
 Av.~Rovisco Pais,
1049-001 Lisboa  (Portugal).}
\email{cordovil@math.ist.utl.pt}
 }
 
\thanks{The  authors' research  was
supported in part by FCT (Portugal) through 
program POCTI and the project  ``Algebraic Methods in  Graph Theory" approved by program Pessoa 2005/6}

\begin{abstract} Let $P$ a  polytope and let $\G(P)$ be the
graph of $P$. Following  Gil Kalai,   we say that
an acyclic orientation $\Or$ of $\G(P)$  is {\bf good} if, for every 
non-empty face $F$ of $P$, the induced graph $\G(F)$ has exactly
one sink.  Gil Kalai  gave a simple way to tell a simple polytope from the good   orientations of its graph. This article is a broader study
of ``good orientations'' (of the graphs) on matroid polytopes.
 \end{abstract}
\maketitle
\centerline{\it Dedicated to Michel Las Vergnas on the occasion of his 65th birthday}

\section{Introduction}
 Let $P$ be a simple polytope (see \cite{Ziegler} for details on polytopes) and let $\G(P)$ be its graph (1-skeleton).
M. Perles conjectured (see the reference [P] in \cite{kalai}) and Blind and Mani \cite{blind}  proved that the graph $\G(P)$ determines the lattice of faces of  $P$.   Kalai \cite{kalai} gave a short and constructive prove of this result; see also \cite{Kaibel, Kaibel2} for a  discussion and refinement of Kalai's technique.  Kalai's proof is based on an intrinsical characterization of the ``good'' orientations of $\G(P)$ between all the acyclic orientations.
 Following \cite{kalai},   we say that
an acyclic orientation $\Or$ of $\G(P)$  is {\bf good} if for every 
non-empty face $F$ of $P$ the induced graph $\G(F)$ has exactly
one {\bf sink}, i.e.,  a vertex of  $\G(F)$ with no lower adjacent vertices.
 Every linear ordering, $\{v_{1}\prec \cdots \prec v_{n}\}$, of the vertex set $V$ of $\G(\M)$
  induces an acyclic orientation $\Or_{\prec}$ of the graph,
where an edge is directed from its larger end-node to its  smaller end-node.  The linear ordering $\{v_{1}\prec \cdots \prec v_{n}\}$ is called  {\bf good}
if  $\Or_{\prec}$ is a good  orientation.
Each acyclic orientation  of an arbitrary graph $G$ is induced by some linear ordering of its vertices, see \cite[Proposition 1.2]{klee2}.
 Good orderings  are in 1-1 correspondence  with 
shelling orderings  of the facets of  the boundary $\partial P^{\Delta}$ of
the dual polytope $P^{\Delta}$ (see Theorem~\ref{shelling2} below for a  matroidal generalization). \par
We say that an   oriented matroid $\M$  is a {\bf  matroid polytope}
if it is acyclic and all the elements of the ground set  $E(\M)$ are extreme points of $\M$. 
 The {\bf graph} $\G(\M)$ of the matroid polytope $\M$ is  the graph whose vertices [resp. edges] are the faces of rank 1 [resp. 2] of $\M$. In particular the vertex set of  $\G(\M)$ is the ground set $E(\M)$.
We say that the   matroid polytope $\M$ is {\bf simple} if every vertex of $\G(\M)$  
is incident with exactly $\text{rank}(\M)-1$ edges. 
The proof of Gil Kalai  \cite{kalai} remains true if we replace ``simple polytope"  by
``simple matroid polytope" so,  the graph
$\G(\M)$ also
encloses the face lattice of $\M$.
If $F$ is a non-empty face of $\M$ then $F$  also is a 
 matroid polytope (we are identifying  $F$ and the restriction of $\M$ to $F$).
For more details on oriented matroid theory see  Section 2 and \cite{OM}.
 \section{K-orderings}
  Let $\M^\prime, \M^{\prime\prime}$ [resp. $\M$] denotes 
  an acyclic oriented matroid [resp. matroid polytope] of rank $r$ on   ground set $E=\{e_1,\dotsc, e_n\}$.  Let $\Or=\Or(\M')$ be
  the class of all the acyclic reorientations of $\M'$.
   For every linear ordering $E(\M)_{\prec}$ let $\Or_{\prec}$ denotes the orientation of $\G(\M)$ where $ \overrightarrow{uv}$ is a directed edge  from $u$ to $v$, if  $v\prec u$.  Let $\G_{\prec}(\M):=(\G(\M), \Or_{\prec})$ be the corresponding digraph. If $F$ is a face of $\M$
let $\G_{\prec}(F)$ be the induced directed subgraph on $F$.
 Note that an element  $e$  is the unique sink of $\G_\prec(\M)$ if and only if for every element $e'$ there 
 is  a directed path, $e'\rightsquigarrow e$, from $e'$ to $e$.\par  
 \begin{definition}\label{shelling ordering}
 {We say that the  linear  ordering 
  $\{e_1\prec  e_2\prec  \cdots \prec e_n\}$ of the ground set of a matroid polytope
  $\M$
   is  a {\bf K-ordering}   if, for every   non-empty face $F$ of $\M$, the directed subgraph $\G_{\prec}(F)$ has exactly 
one sink. In particular the element $e_1$ is the unique sink of the digraph $\G_{\prec}(\M)$.
}
\end{definition}
 From  the  ``Topological Representation Theorem",  we know that there is a pure regular CW-complex  of dimension $r-1$,  
  with the topology of a PL-sphere and encoding $\Or$, see \cite[Theorem 5.2.1]{OM}. We will denote this CW-complex by $\Delta(\Or)$, or  by $\Delta$ for short, and called it   the
  {\bf $CW$-complex of   acyclic reorientations of $\M'$}. 
 To  every cell $W\in \Delta$ we
  attach   a ``sign vector''   $\sigma(W)\in \{0,+,-\}^{E(\M')}$ (for precisions see \cite{OM}). We will identify $W$ and $\sigma(W)$ and  set $W_{(e_i)}=\sigma(W)_{(e_i)}$.
 The set $$\text{supp}(W)=\{e_i:  W_{(e_i)}\not =0\}$$ is called  the  {\bf support} of  $W$. We say that 
 the cell $W'$ is a face of  $W$ if the following two conditions hold:
  \begin{itemize}\label{face}
  \item[$(\ref{face}.1)$]
  $\text{supp}(W')\subset \text{supp}(W)$;
  \item[$(\ref{face}.2)$] For every $e_i\in E$, if  $e_i\in \text{supp}(W')$ then we have  $W'_{(e_i)}=W_{(e_i)}$.   
  \end{itemize}
 A cell  $P\in \Delta$  has dimension 0 (is a {\bf vertex}) if, seen as a signed vector, $P$ 
is a signed cocircuit of $\M$.    We will see $\Delta$ as an ``abstract'' regular
 cell complex over the set of its vertices. Every cell $W\in \Delta$ is identified  with the set of its vertices $V(W)$:
  $$W\equiv V(W):=\{P:\, P\,\,\, \text{is a vertex of}\,\, \Delta\,\, \text{and}\,\, P \leq W\}.
 $$
 The facets   of $\Delta$
    are called the  {\bf topes}  and have   support equal to $E(\M)$.
   We use the letter $T$ to denote a tope. Every tope   $T\in \Delta$ is a  $PL$-ball and its boundary $\partial T$ a $PL$-sphere.
  Every tope $T$ fixes  one acyclic reorientation $\M^{\prime\prime}\in \Or(\M')$. To explicit this correspondence we write $T=T(\M^{\prime\prime})$. (The opposite tope $-T$ fix the same oriented matroid $\M^{\prime\prime}$.) If a fixed  acyclic oriented matroid $\M'$ is given    we set $T(\M')=(+,+,\dotsc, +)$.   Let $\La(\M)$  [resp. $\La(T(\M))$] be 
   the lattice of the faces (with both trivial faces) of $\M$ [resp. $T(\M)$]  ordered by inclusion.  There is a canonical
    anti-isomorphism  ${\Xi}: \La(\M )\to\La(T(\M))$
where, to every element $X$ of rank $s$ of $\La(\M)$,   $\Xi(X)$ denotes the element of rank $r-s$ of $\La (T(\M))$ determined by the following conditions:
 $$
 \Xi(X)_{(e_\ell)}= \begin{cases}
 0 & \quad \text{if}  \quad e_\ell\leq X\\
 + & \quad \text{if}  \quad e_\ell\not \leq X.
 \end{cases}
  $$
  The atoms of $\La(T(\M))$ are the image by $\Xi$ of the co-atoms
  of $\La(\M)$. The set of vertices and facets of   CW-complex $T$, are respectively
  $$\{\Xi(H):\; H\, \text{a facet of}\,\, \M\}\quad\text{and}\quad\{\Xi (e):\; e\in E(\M)\}.$$
If $\M$ is a simple matroid polytope  then $T(\M)$
is an (abstract) simplicial complex of dimension $\text{rank}(\M)-1$.  In particular,
for every pair $\{W, W'\}$ of cells of $T(\M)$, we have
$V(W\wedge W')= V(W)\cap V(W')$.
The following definition is a particular case of the standard one, see \cite{klee,kleins, Ziegler} for details.
(The equivalence of   Conditions $(\ref{sh}.1)$
 and $(\ref{sh}.1')$ is  left to the reader.)
     \begin{definition}\label{sh}
 Let $ \Delta$  be
the $CW$-complex of the acyclic reorientations  of $\M$ and  $T\in \Delta$ be the tope associated to $\M$.
 We say that    the linear ordering   $\{\Xi_1:=\Xi(e_1)\prec \Xi_2\prec \cdots \prec \Xi_n\}$  is a {\bf shelling} of the $PL$-sphere $\partial T$ and  $\partial T$ is {\bf shellable} if one of the following  equivalent conditions holds:
 \begin{itemize}
 \item[$(\ref{sh}.1)$] For every pair of co-atoms $\Xi_i\prec \Xi_j$ such that $\Xi_j\cap \Xi_i\not=\emptyset$, there is some facet  $\Xi_\ell\prec \Xi_j$ such that $\Xi_\ell\cap \Xi_j$ is an abstract  simplex of cardinality $r-2$ of $T$
 and $\Xi_j\cap \Xi_i\subseteq \Xi_\ell\cap \Xi_j$;
  \item[$(\ref{sh}.1')$] For every pair of vertices $\{e_i, e_j\}$, $1\leq i<j\leq n$, on a non-singular face $F$ of $\M$, there is some  $\ell<j$ such that $ \overrightarrow{e_je_\ell}$
  is a directed edge of the digraph $\G_\prec(F)$.
 \end{itemize}
 \end{definition}
  \begin{theorem}\label{shelling2}
  The  following two statements are equivalent:
\begin{itemize}
\item[$(\ref{shelling2}.1)$]
$\{e_1\prec  \cdots \prec e_n\}$ is a
 K-ordering  of $\M$;
\item[$(\ref{shelling2}.2)$]
$\{ \Xi_1=    \Xi(e_1)\prec  \Xi_2\prec \cdots \prec \Xi_n\} $
 is a shelling of the $PL$-sphere $\partial T(\M)$.
  \end{itemize}
 \end{theorem}
  \begin{proof}
   \noindent$(\ref{shelling2}.1)\Longrightarrow (\ref{shelling2}.2)$.
  As $i<j$ we know that $e_j$  is not a sink of $\G_\prec(F)$. So there is a directed edge $ \overrightarrow{e_je_\ell}$ of $\G_\prec(F)$
   and $(\ref{sh}.1')$ holds.\par
   \noindent $(\ref{shelling2}.2)\Longrightarrow (\ref{shelling2}.1)$.
Let
  $\{e_{i_1}\prec\cdots\prec e_{i_{n'}}\}$
  be the induced ordering on a non-empty face $F$  of $\M$. From Condition~$(\ref{sh}.1')$
  we know that $e_{i_1}$ is the unique sink of $\G_\prec(F)$ and $(\ref{shelling2}.1)$ follows.
 \end{proof}
 Let us recall that a linear ordering  $\{e_1\prec  e_2\prec  \cdots \prec e_n\}$ of the  ground set $E(\M')$ is called a {\bf shelling ordering}  of  $\M'$, if the orientation obtained from $\M'$  by  changing the signs on the initial sets $E_{i}:=\{e_1,e_2,\dotsc, e_i\}$, $i=1, 2,\dotsc, n$,  is also  acyclic. Edmonds and Mandel proved that in this case,
 $\{ \Xi_1=      \Xi(e_1)\prec  \Xi_2\prec \cdots \prec \Xi_n\} $
 is a shelling of the  PL-sphere $\partial T(\M)$, see \cite[Proposition~4.3.1]{OM}.
Note that if $F$ is a face of a matroid polytope $\M$, every shelling ordering of $\M$ induces a shelling ordering on $F$.
 The following result is a consequence of Theorem~\ref{shelling2} and the above result  of Edmonds and Mandel.
 \begin{corollary}
 Every  shelling ordering of a simple matroid polytope is also a K-ordering.
 \qed
 \end{corollary}
 Let $f_\ell(\M)$ be the number of faces of rank $\ell+1$ of $\M$, $-1\leq \ell \leq r-1$. By convention set 
 $f_{-1}(\M)=f_{r-1}(\M)=1$.
 By analogy with the definition of  the $\mathbf{f}$-vector and the $\mathbf{h}$-vector of a polytope,  we call  the  vector,
 $$\mathbf{f}(\M):=\;\big(f_{-1}(\M), f_0(\M), f_1(\M),\dotsc, f_{r-2}(\M), f_{r-1}(\M)\big),$$ the $\mathbf{f}$-{\bf vector} of the  (simple) matroid  polytope $\M$  and we call  the vector,  $$\mathbf{h}^*(\M):=\;\big(h^*_0(\M),h^*_1(\M),\dotsc, h^*_{r-1}(\M)\big),$$ 
 determined by  the  formulas
 \begin{equation}\label{h**}
 {h}^*_\ell(\M)=\sum_{i=0}^\ell(-1)^{\ell-i}\binom{r-1-i}{\ell- i} f_{r-1-i}(\M),\,\,\, \ell=0, 1,\dotsc, r-1,
 \end{equation}
  the $\mathbf{h}^*$-{\bf vector} of $\M$.
Note that  the $\mathbf{f}$-vector also can be recovered from the  $\mathbf{h}^*$-vector:
\begin{equation}\label{h*}
 f_\ell(\M)=\sum_{i=0}^{r-1-\ell} \binom{r-1-i}{\ell} h_i^*(\M),\,\,\, \ell=0,1,\dotsc, r-1.
\end{equation}
 (See \cite[Section 8.3]{Ziegler} for a good survey of $\mathbf{h}$-vectors of simplicial polytopes and Dehn-Sommerville Equations.)
 From   Euler-Poincar\'e formula (see   
\cite[Corollary 4.6.11]{OM}) we know that
\begin{align}\label{Euler}
\sum_{i=-1}^{r-1} (-1)^i f_i(\M)=0.
\end{align}
The graph $\G(\M)$ is regular of degree $r-1$. Fix a K-ordering $\{e_1\prec  \cdots \prec e_n\}$ of   $\M$. Let $d^+ (e)$
 [resp. $d^{-} (e)$] denotes the outdegree [resp. indegree] of $e\in E$.  Set 
 $$\mathbf{d}_\ell^+(\M):=\; |\{e: e\in E,\, d^+ (e)=\ell\}|,$$
  $$\mathbf{d}_\ell^-(\M) :=\; |\{e: e\in E,\, d^- (e)=\ell\}|.$$
We clearly have that $\mathbf{d}_\ell^+(\M)= \mathbf{d}_{r-1-\ell}^-(\M)$.
 \begin{theorem}\label{final} 
 The integers $\mathbf{d}_\ell^+(\M)$, $\ell=0, 1,\dotsc, r-1$, are invariant of 
 the matroid polytope $\M$
 (i.e., are independent of the K-ordering)  and they are determined by the equalities:
 \begin{equation}\label{essenc}
  h_\ell^*(\M)=\mathbf{d}_\ell^+(\M), \,\, \ell=0, 1,\dotsc, r-1.
\end{equation}  
\end{theorem}
 \begin{proof}  Note that 
\begin{align*}
\mathbf{f}_\ell(\M)
  &=\sum_{j=\ell}^{r-1} \binom{j}{\ell}\mathbf{d}^-_j (\M)\\
 &= \sum_{i=0}^{r-1-\ell}  \binom{r-1-i}{\ell}\mathbf{d}^-_{r-1-i}(\M)\\
  &= \sum_{i=0}^{r-1-\ell} \binom{r-1-i}{\ell}\mathbf{d}^+_{i} (\M).
 \end{align*}
  From Equation~$(\ref{h*})$ we conclude that $ \mathbf{d}^+_{\ell}(\M)={h}^*_\ell(\M)$.
 \end{proof}
 \begin{corollary}\label{reverse} 
Let  $E_\prec:=\{e_1\prec  e_2\prec  \cdots \prec e_n\}$
be a K-ordering of $\M$.
 Then the  reverse ordering 
 $E_{\prec^*}:=\{e_n\prec^*  e_{(n-1)}\prec^*  \cdots \prec^* e_1\}$
 is also a K-ordering of  $\M$.
 \end{corollary}
 \begin{proof} 
 It is necessary to prove that, for every non-empty face $F$ of $\M$,  the digraph 
 $\G_{\prec^*}(F)$ has exactly one  sink, i.e., $\mathbf{d}_{\text{rank}(F)-1}^+(F)=1$.
 From the equalities $(\ref{essenc})$  we know $\mathbf{d}_{r-1}^+(\M)=h^*_{r-1}(\M)$.
    From  Euler-Poincar\'e formula $(\ref{Euler})$ we conclude  that  
 \begin{align} h^*_{r-1}(\M) & = \sum_{i=0}^{r-1} (-1)^{r-1-i}  f_{r-1-i}(\M)=\\
  & =  \sum_{j=0}^{r-1} (-1)^{j}  f_{j}(\M)= f_{-1}(\M)=1.
 \end{align}
So 
   $\G_{\prec^*}(M)$ has exactly  one sink. 
  As every face $F$ of $\M$
   is a simple matroid polytope the result follows.
  \end{proof}
\section{The Cube}
  Let $C^{d}:=\; \{x\in \R^{d}:\; 0\leq x_{\ell} \leq 1,\,\, \ell=1, \dotsc, d\}$ 
be the $d$-dimensional cube.
 As the polar of the  cube $C^{d}$ is the $d$-dimensional crosspolytope
 (a simplicial polytope), it results from Theorem~\ref{final} above that:
  \begin{equation}
 \mathbf{h}^*(C^d)=\Bigg( \binom{d}{0}, \binom{d}{1},\dotsc,\binom{d}
  {d-1},\binom{d}{d}\Bigg),
   \end{equation}
   i.e., there are exactly $\binom{d}{\ell}$, $ 0\leq \ell\leq d$, vertices of $\G(C^d)$ such that $\mathbf{d}^-(e)=d-\ell$.\par
 Let $B:=\{0, 1\}^d$ be the set  of the vertices of the  cube $C^{d}$ and for every
$b_i\in B$, $1\leq i\leq 2^d$, set $e_i=(b_i, 1)\in \R^{d+1}$. The rank $d+1$  {\bf cube matroid polytope}, $\C^{d}$, is the oriented matroid  determined by
 the linear dependencies of   vectors of $E:= \{e_{i}: b_{i}\in B\}$.
  The following theorem is closely related to the results presented here. (We present here a  slightly different version of the original result.)
  \begin{theorem} \cite[Proposition 3.2]{hammer}\label{known}
  Let $\C^d$ be the cube matroid polytope of rank at least three.  Let   $\{e_1\prec  e_2\prec  \cdots \prec e_{2^d}\}$
be a linear ordering of $E(\C^d)$. Then the following two conditions are equivalent:
\begin{itemize}
\item[$(\ref{known}.1)$]
$\{e_1\prec  e_2\prec  \cdots \prec e_{2^d}\}$
 is a K-ordering of\, $\C^d$;
\item[$(\ref{known}.2)$] For every rank three face $F$ of\, $\C^d$, 
the  digraph $\G_\prec(F)$ has an unique sink.
\qed
\end{itemize}
\end{theorem}
In the  rank three cube  matroid polytope $\C^2$, the ``K-ordering'' and the ``shelling orderings''  coincide.   
This result suggest the following problem:
\begin{problem}  Is there  a
simple characterisation of  shelling orderings of the   cube matroid polytope\, $\C^d$?
Are there  K-orderings of the cube matroid polytope that  are not shelling orderings? 
\end{problem}

\end{document}